\documentclass[12pt]{article}%
\usepackage{graphicx}
\usepackage[intlimits]{amsmath}
\usepackage{latexsym}
\usepackage{amsfonts}
\usepackage{amssymb}%
\makeatletter
\newfont{\footsc}{cmcsc10 at 8truept}
\newfont{\footbf}{cmbx10 at 8truept}
\newfont{\footrm}{cmr10 at 10truept}
\newtheorem{theorem}{Theorem}

\newtheorem{proposition}[theorem]{Proposition}

\newenvironment{proof}[1][Proof]{\noindent{\textbf {#1}  }}  {\hfill$\Box$\bigskip}

\begin{document}

\title{Regular, pseudo-regular, and almost regular matrices}
\author{Vladimir Nikiforov\\{\small Department of Mathematical Sciences, University of Memphis, }\\{\small Memphis TN 38152, USA, email: }\textit{vnikifrv@memphis.edu}}
\maketitle

\begin{abstract}
We give lower bounds on the largest singular value of arbitrary matrices, some
of which are asymptotically tight for almost all matrices. To study when these
bounds are exact, we introduce several combinatorial concepts. In particular,
we introduce regular, pseudo-regular, and almost regular matrices.
Nonnegative, symmetric, almost regular matrices were studied earlier by
Hoffman, Wolfe, and Hoffmeister.\medskip


\end{abstract}

\section{Introduction}

Let $\Sigma\left(  A\right)  $ be the sum of the entries of a matrix $A.$
Hoffman, Wolfe and Hofmeister \cite{HWH95} showed that if $A=\left(
a_{ij}\right)  $ is a nonnegative symmetric matrix with positive rowsums
$d_{1},\ldots,d_{n}$ and $\mu\left(  A\right)  $ is its largest eigenvalue,
then
\begin{equation}
\mu\left(  A\right)  \geq\frac{1}{\Sigma\left(  A\right)  }%
{\textstyle\sum\limits_{i,j}}
a_{ij}\sqrt{d_{i}d_{j}}\label{hwhin}%
\end{equation}
with equality holding if and only if $d_{i}d_{j}=\mu^{2}\left(  A\right)  $
whenever $a_{ij}>0.$

The aim of this note is to extend this result in several directions. First,
instead of $\mu\left(  A\right)  ,$ we consider the largest singular value
$\sigma\left(  A\right)  ,$ thereby dropping the requirement that $A$ is
symmetric, square, and nonnegative. Second, we present wider classes of lower
bounds on $\sigma\left(  A\right)  $ some of which are asymptotically tight
for almost all matrices. Finally, we study when these lower bounds are exact,
thus introducing regular, pseudo-regular, and almost regular matrices. We also
introduce a few combinatorial concepts to support the study of structural
properties of arbitrary matrices.

For basic notation and definitions see \cite{HoJo88}. Specifically, we call a
matrix \emph{scalar} if it is a scalar multiple of a nonnegative matrix. Also,
we write $\mathbf{j}_{m}$ for the vector of $m$ ones.

\section{Main results}

Let $A=\left(  a_{ij}\right)  $ be an $m\times n$ matrix with row and column
sums $r_{1},\ldots,r_{m}$ and $c_{1},\ldots,c_{n}.$ We first generalize the
values $c_{i}$ and $r_{j}.$ Index the rows and columns of $A$ by the elements
of two disjoint sets $R=R\left(  A\right)  $ and $C=C\left(  A\right)  .$ For
all $i\in R\cup C,$ set $w_{A}^{1}\left(  i\right)  =1;$ for all $s\geq2,$
$i\in R,$ $j\in C,$ set
\[
w_{A}^{s}\left(  i\right)  =%
{\textstyle\sum\limits_{k\in C}}
a_{ik}w_{A}^{s-1}\left(  k\right)  ,\text{ \ \ \ \ }w_{A}^{s}\left(  j\right)
=%
{\textstyle\sum\limits_{k\in R}}
a_{kj}w_{A}^{s-1}\left(  k\right)  .
\]
Finally, for all $s\geq1,$ set
\[
w_{A}^{s}\left(  R\right)  =%
{\textstyle\sum\limits_{k\in R}}
w_{A}^{s}\left(  k\right)  ,\text{ \ \ \ \ }w_{A}^{s}\left(  C\right)  =%
{\textstyle\sum\limits_{k\in C}}
w_{A}^{s}\left(  k\right)  .
\]

Note that $w_{A}^{2}\left(  i\right)  =r_{i}$ if $i\in R,$ and $w_{A}%
^{2}\left(  i\right)  =c_{i}$ if $i\in C$. Also, if $A$ is the adjacency
matrix of a graph, $w_{A}^{s}\left(  i\right)  $ is the number of walks on $s$
vertices starting with the vertex $i.$

Using somewhat different notation, in \cite{Nik07} it is proved that for every
$m\times n$ matrix $A$ and all odd $p$ and $r$ such that $p>r\geq1,$%
\begin{equation}
\sigma^{p-r}\left(  A\right)  w_{A}^{r}\left(  R\right)  \geq w_{A}^{p}\left(
R\right)  . \label{mainin}%
\end{equation}
Moreover, for all $s\geq1,$
\begin{equation}
\sigma^{2s}\left(  A\right)  =\lim_{r\rightarrow\infty}\frac{w_{A}%
^{2r+2s+1}\left(  R\right)  }{w_{A}^{2r+1}\left(  R\right)  }=\lim
_{r\rightarrow\infty}\max_{k\in R\left(  A\right)  }\frac{w_{A}^{2r+2s+1}%
\left(  k\right)  }{w_{A}^{2r+1}\left(  k\right)  } \label{maineq}%
\end{equation}
unless the eigenspace of $AA^{\ast}$ corresponding to $\sigma^{2}\left(
A\right)  $ is orthogonal to $\mathbf{j}_{m}.$

Note that: \emph{(i)} inequality (\ref{mainin}) may not hold if $p$ or $r$ are
even (see \cite{Bru06}, p. 728\ and \cite{Nik06}, p. 262); \emph{(ii)}
equalities (\ref{maineq}) hold if $A$ is a nonzero scalar matrix; \emph{(iii)}
inequality (\ref{mainin}) implies a number of known results on the spectral
radius of graphs (see \cite{Nik06}, p. 258).

Inequality (\ref{mainin}) can be proved using the Rayleigh principle. This
simple approach helps produce other similar bounds of increasing complexity.
We shall focus on the following general inequality.

\begin{theorem}
\label{th1}Let $A$ be a matrix, $R=R\left(  A\right)  $ and $C=C\left(
A\right)  .$ Then for all $r\geq1,$%
\begin{equation}
\sigma\left(  A\right)  \sqrt{w_{\left\vert A\right\vert }^{r}\left(
R\right)  w_{\left\vert A\right\vert }^{r}\left(  C\right)  }\geq\left\vert
{\textstyle\sum\limits_{i\in R,j\in C}}
a_{ij}\sqrt{w_{\left\vert A\right\vert }^{r}\left(  i\right)  w_{\left\vert
A\right\vert }^{r}\left(  j\right)  }\right\vert . \label{secin}%
\end{equation}

\end{theorem}

Particularly, for $r=1$ Theorem \ref{th1} reads as
\begin{equation}
\sigma\left(  A\right)  \geq\left\vert \Sigma\left(  A\right)  \right\vert
/\sqrt{nm}. \label{secin1}%
\end{equation}
Also, since $w_{\left\vert A\right\vert }^{2}\left(  R\right)  =w_{\left\vert
A\right\vert }^{2}\left(  C\right)  =\Sigma\left(  \left\vert A\right\vert
\right)  $, for $r=2$ Theorem \ref{th1} extends inequality (\ref{hwhin}) to
\[
\sigma\left(  A\right)  \Sigma\left(  \left\vert A\right\vert \right)
\geq\left\vert
{\textstyle\sum\limits_{i\in R,j\in C}}
a_{ij}\sqrt{w_{\left\vert A\right\vert }^{2}\left(  i\right)  w_{\left\vert
A\right\vert }^{2}\left(  j\right)  }\right\vert .
\]

It is natural to study when equality holds in inequalities (\ref{mainin}) and
(\ref{secin}). To this end we first introduce some combinatorial concepts.

\subsection{A few combinatorial concepts}

For any matrix $A=\left(  a_{ij}\right)  ,$ let $\mathcal{B}\left(  A\right)
$ be the bipartite graph with vertex classes $R\left(  A\right)  $ and
$C\left(  A\right)  $ such that $i\in R\left(  A\right)  $ is joined to $j\in
C\left(  A\right)  $ whenever $a_{ij}\neq0.$

Call a matrix $A$ \emph{connected} if $\mathcal{B}\left(  A\right)  $ is
connected. Note that a symmetric matrix is connected exactly when it is irreducible.

For scalar matrices connectedness can be expressed in terms of their powers.

\begin{proposition}
\label{pro0}A scalar matrix $A$ is connected if and only if for every $i\in
R\left(  A\right)  ,$ $j\in C\left(  A\right)  ,$ there exists $r$ such that
the $\left(  i,j\right)  $ entry of $\left(  AA^{\ast}\right)  ^{r}A$ is nonzero.
\end{proposition}

Call a maximal connected submatrix of $A$ a \emph{component} of $A$.

We say that $A$ is \emph{cogredient} to $B$ if there exist permutation
matrices $P$ and $Q$ such that $A=PBQ.$

The following two assertions are obvious.

\begin{proposition}
\label{pro1}If a matrix $A$ has no zero rows or columns, then it is cogredient
to a block diagonal matrix%
\[
\left(
\begin{array}
[c]{ccc}%
A_{1} & \cdots & 0\\
\vdots & \ddots & \vdots\\
0 & \cdots & A_{r}%
\end{array}
\right)
\]
where $A_{1},\ldots,A_{r}$ are the components of $A.$
\end{proposition}

\begin{proposition}
\label{pro2}The multiset of the nonzero singular values of $A$ is the union of
the multisets of the nonzero singular values of its components. In particular,
for every matrix $A,$%
\[
\sigma\left(  A\right)  =\max\left\{  \sigma\left(  C\right)  :C\text{ is a
component of }A\right\}  .
\]

\end{proposition}

Let $A$ be a nonzero scalar matrix. We call $A$ \emph{regular} if its row sums
are equal and so are its columns sums.

We call $A$ \emph{pseudo-regular }if\emph{ }$w_{A}^{5}\left(  i\right)
=\lambda w_{A}^{3}\left(  i\right)  $ for all $i\in R\left(  A\right)  $ and
fixed $\lambda.$ Equivalently, $A$ is pseudo-regular if the vector with
coordinates $w_{A}^{3}\left(  i\right)  ,$ $i\in R\left(  A\right)  $ is an
eigenvector of $AA^{\ast}$.

If each component $C$ of $A$ is regular and $\sigma\left(  C\right)
=\sigma\left(  A\right)  ,$ we call $A$ \emph{almost-regular. }

Note that regular matrices generalize doubly stochastic matrices.
Pseudo-regular matrices generalize the adjacency matrices of pseudo-regular
and pseudo-semiregular graphs (see \cite{Row07} for a comprehensive survey).
Almost regular matrices extend the concept introduced for nonnegative
symmetric matrices in \cite{HWH95}.

It is easy to see that regular matrices are almost regular, and that almost
regular matrices are pseudo-regular. However the matrix%
\[
A=\left(
\begin{array}
[c]{cccc}%
1 & 1 & 0 & 0\\
1 & 0 & 1 & 0\\
1 & 0 & 0 & 1
\end{array}
\right)
\]
is connected and pseudo-regular, but not regular. Note also that $A^{\ast}$ is
not pseudo-regular.

Here is a complete characterization of pseudo-regular matrices.

\begin{proposition}
\label{pro3} A nonzero $m\times n$ scalar matrix $A$ is pseudo-regular if and
only if the following conditions hold:

(i) the vector with entries $w_{A}^{3}\left(  i\right)  ,$ $i\in R\left(
A\right)  $ is an eigenvector of $AA^{\ast}$ to some nonzero eigenvalue
$\mu\left(  AA^{\ast}\right)  ;$

(ii) the eigenvectors of $AA^{\ast}$ to every nonzero eigenvalue $\mu^{\prime
}\left(  AA^{\ast}\right)  \neq\mu\left(  AA^{\ast}\right)  $ are orthogonal
to $\mathbf{j}_{m}.$
\end{proposition}

Using this characterization, we can relax the definition of pseudo-regularity,
preserving the same property scope.

\begin{proposition}
\label{pro4}Suppose that $A$ is a scalar matrix, $r,s$ are odd, and
$r>s\geq3.$ If $w_{A}^{r}\left(  i\right)  =\lambda w_{A}^{s}\left(  i\right)
$ for all $i\in R\left(  A\right)  $ and fixed $\lambda,$ then $A$ is pseudo-regular.
\end{proposition}

\bigskip

\subsection{Sufficient conditions for equality in (\ref{mainin}) and
(\ref{secin})}

The following theorem gives a condition for equality in (\ref{mainin}).

\begin{theorem}
\label{th2}Suppose that $A$ is a scalar matrix with $R=R\left(  A\right)  .$
If
\[
\sigma^{2s}\left(  A\right)  w_{A}^{2r+1}\left(  R\right)  =w_{A}%
^{2r+2s+1}\left(  R\right)  .
\]
for some $s\geq1,$ $r\geq0,$ then $A$ is pseudo-regular.
\end{theorem}

Similar double condition implies a stronger conclusion.

\begin{theorem}
\label{th2.1}Suppose that $A$ is a scalar matrix with $R=R\left(  A\right)  ,$
$C=C\left(  A\right)  .$ If
\begin{align}
\sigma^{2s}\left(  A\right)  w_{A}^{1}\left(  R\right)   &  =w_{A}%
^{2s+1}\left(  R\right)  ,\label{eqr}\\
\sigma^{2r}\left(  A\right)  w_{A}^{1}\left(  C\right)   &  =w_{A}%
^{2r+1}\left(  C\right)  \label{eqc}%
\end{align}
for some $r,s\geq1,$ then $A$ is almost regular.
\end{theorem}

Next, we generalize the second part of the aforementioned theorem of Hoffman,
Wolfe, and Hofmeister giving conditions for equality in (\ref{secin}).

\begin{theorem}
\label{th3}Let $A=\left(  a_{ij}\right)  $ be a scalar matrix and $r\geq1,$
$s\geq1$. The following three conditions are equivalent:

(i) $A$ is almost regular;

(ii) $\left\vert w_{A}^{r}\left(  i\right)  w_{A}^{r}\left(  j\right)
\right\vert =\sigma\left(  A\right)  ^{2}\left\vert w_{A}^{r}\left(  R\right)
w_{A}^{r}\left(  C\right)  \right\vert $ whenever $a_{ij}\neq0;$

(iii) we have
\begin{equation}
\sigma\left(  A\right)  \sqrt{\left\vert w_{A}^{r}\left(  R\right)  w_{A}%
^{r}\left(  C\right)  \right\vert }=\left\vert
{\textstyle\sum\limits_{i\in R,j\in C}}
a_{ij}\sqrt{\left\vert w_{A}^{r}\left(  i\right)  w_{A}^{r}\left(  j\right)
\right\vert }\right\vert .\label{eq}%
\end{equation}

\end{theorem}

A stronger condition holds for equality in (\ref{secin}) with $r=1$.

\begin{theorem}
\label{th4}A scalar matrix $A\in M_{m,n}$ is regular if and only if
$\sigma\left(  A\right)  =\left\vert \Sigma\left(  A\right)  \right\vert
/\sqrt{nm}$.
\end{theorem}

Note that the assumption that $A$ is scalar is essential in Theorems \ref{th2}
to \ref{th4}. Indeed letting
\[
A=\left(
\begin{array}
[c]{cc}%
1+i & 1-i\\
1-i & 1+i
\end{array}
\right)  ,
\]
we see that
\begin{align*}
\sigma\left(  A\right)   &  =2,\text{ \ \ }\Sigma\left(  A\right)  =4,\text{
\ \ }w_{A}^{1}\left(  R\right)  =w_{A}^{1}\left(  C\right)  =2,\\
w_{A}^{2}\left(  R\right)   &  =w_{A}^{2}\left(  C\right)  =4,\text{
\ \ }w_{A}^{3}\left(  R\right)  =8,\\%
{\textstyle\sum\limits_{i\in R,j\in C}}
a_{ij}\sqrt{\left\vert w_{A}^{2}\left(  i\right)  w_{A}^{2}\left(  j\right)
\right\vert }  &  =%
{\textstyle\sum\limits_{i\in R,j\in C}}
2a_{ij}=8.
\end{align*}
Thus, we have
\begin{align*}
\sigma^{2}\left(  A\right)  w_{A}^{1}\left(  R\right)   &  =w_{A}^{3}\left(
R\right)  ,\\
\sigma\left(  A\right)  \sqrt{\left\vert w_{A}^{2}\left(  R\right)  w_{A}%
^{2}\left(  C\right)  \right\vert }  &  =\left\vert
{\textstyle\sum\limits_{i\in R,j\in C}}
a_{ij}\sqrt{\left\vert w_{A}^{2}\left(  i\right)  w_{A}^{2}\left(  j\right)
\right\vert }\right\vert ,\\
\sigma\left(  A\right)   &  =\left\vert \Sigma\left(  A\right)  \right\vert
/2,
\end{align*}
although $A$ is not scalar.

\section{\label{prf}Proofs}

\begin{proof}
[\textbf{Proof of Theorem \ref{th1}}]Set $x_{i}=\sqrt{w_{\left\vert
A\right\vert }^{r}\left(  i\right)  /w_{\left\vert A\right\vert }^{r}\left(
R\right)  }$ for all $i\in R$ and let $\mathbf{x}=\left(  x_{i}\right)  .$
Likewise, set $y_{i}=\sqrt{w_{\left\vert A\right\vert }^{r}\left(  i\right)
/w_{\left\vert A\right\vert }^{r}\left(  C\right)  }$ for all $i\in C$ and let
$\mathbf{y}=\left(  y_{i}\right)  .$ Since $\left\Vert \mathbf{x}\right\Vert
=\left\Vert \mathbf{y}\right\Vert =1,$ by Schur's lemma \cite{Sch11}, we
obtain%
\begin{align*}
\sigma\left(  A\right)    & =\max_{\left\Vert \mathbf{u}\right\Vert
=\left\Vert \mathbf{v}\right\Vert =1}\left\vert \left\langle A\mathbf{u}%
,\mathbf{v}\right\rangle \right\vert \geq\left\vert \left\langle
A\mathbf{x},\mathbf{y}\right\rangle \right\vert \\
& =\frac{1}{\sqrt{w_{\left\vert A\right\vert }^{r}\left(  R\right)
w_{\left\vert A\right\vert }^{r}\left(  C\right)  }}\left\vert
{\textstyle\sum\limits_{i\in R,j\in C}}
a_{ij}\sqrt{w_{\left\vert A\right\vert }^{r}\left(  i\right)  w_{\left\vert
A\right\vert }^{r}\left(  j\right)  }\right\vert ,
\end{align*}
completing the proof.
\end{proof}

\bigskip

In the proofs below we shall assume that $A$ is an $m\times n,$ nonzero,
nonnegative matrix with $R=R\left(  A\right)  $ and $C=C\left(  A\right)  ;$
$r_{1},\ldots,r_{m},$ $c_{1},\ldots,c_{n}$ are its row and column sums, and
$\sigma=\sigma_{1}\geq\cdots\geq\sigma_{m}$ are its singular values.

Let $AA^{\ast}=VDV^{\ast}$ be the unitary decomposition of $AA^{\ast};$ thus,
the columns of $V$ are the unit eigenvectors to $\sigma_{1}^{2},\ldots
,\sigma_{m}^{2}$ and $D$ is the diagonal matrix with $\sigma_{1}^{2}%
,\ldots,\sigma_{m}^{2}$ along its main diagonal. Then for every $l\geq0,$
\[
w_{A}^{2l+1}\left(  R\right)  =\Sigma\left(  \left(  AA^{\ast}\right)
^{l}\right)  =\Sigma\left(  VD^{l}V^{\ast}\right)  =%
{\textstyle\sum\limits_{i\in\left[  m\right]  }}
c_{i}\sigma_{i}^{2l},
\]
where $c_{i}=\left\vert \sum_{j\in\left[  m\right]  }v_{ji}\right\vert
^{2}\geq0$ is independent of $l.$

Note also that for all $r,s\geq0,$%
\begin{equation}%
{\textstyle\sum_{i\in R}}
w_{A}^{2r+1}\left(  i\right)  w_{A}^{2s+1}\left(  i\right)  =w_{A}%
^{2r+2s+1}\left(  R\right)  . \label{ideq}%
\end{equation}
We omit the easy proof by induction on $s$.

\bigskip

\begin{proof}
[\textbf{Prrof of Theorem \ref{th2}}]In the above notation we have
\[
\sigma^{2s}%
{\textstyle\sum\limits_{i\in\left[  m\right]  }}
c_{i}\sigma_{i}^{2r}=\sigma^{2s}w_{A}^{2r+1}\left(  R\right)  =w_{A}%
^{2r+2s+1}\left(  R\right)  =%
{\textstyle\sum\limits_{i\in\left[  m\right]  }}
c_{i}\sigma_{i}^{2r+2s}.
\]
Hence, if $0<\sigma_{i}^{2}<\sigma^{2}$, then $c_{i}=0.$ Therefore, for all
$r>1,$ we have $w_{A}^{2r+1}\left(  R\right)  =C\sigma^{2r},$ where $C$ is
independent of $r.$ Specifically,
\[
w_{A}^{5}\left(  R\right)  =C\sigma^{4},\text{ \ \ }w_{A}^{7}\left(  R\right)
=C\sigma^{6},\text{ \ \ }w_{A}^{9}\left(  R\right)  =C\sigma^{8}.
\]
Note the following instances of identity (\ref{ideq})%
\[
w_{A}^{5}\left(  R\right)  =%
{\textstyle\sum\limits_{k\in R}}
\left(  w_{A}^{3}\left(  k\right)  \right)  ^{2},\text{ \ \ }w_{A}^{7}\left(
R\right)  =%
{\textstyle\sum\limits_{k\in R}}
w_{A}^{5}\left(  k\right)  w_{A}^{3}\left(  k\right)  ,\text{ \ \ }w_{A}%
^{9}\left(  R\right)  =%
{\textstyle\sum\limits_{k\in R}}
\left(  w_{A}^{5}\left(  k\right)  \right)  ^{2}.
\]
Hence, using the Cauchy-Schwarz inequality, we obtain%
\begin{align*}
C\sigma^{6}  &  =w_{A}^{7}\left(  R\right)  =%
{\textstyle\sum\limits_{k\in R}}
w_{A}^{5}\left(  k\right)  w_{A}^{3}\left(  k\right)  \leq\sqrt{%
{\textstyle\sum\limits_{k\in R}}
\left(  w_{A}^{5}\left(  k\right)  \right)  ^{2}%
{\textstyle\sum\limits_{k\in R}}
\left(  w_{A}^{3}\left(  k\right)  \right)  ^{2}}\\
&  =\sqrt{w_{A}^{9}\left(  R\right)  w_{A}^{5}\left(  R\right)  }=C\sigma^{6}.
\end{align*}

We have equality in the Cauchy-Schwarz inequality; hence for each $k\in R,$
$w_{A}^{5}\left(  k\right)  =\lambda w_{A}^{3}\left(  k\right)  ,$ where
$\lambda$ is independent of $k.$ Therefore $A$ is pseudo-regular, completing
the proof.
\end{proof}

\bigskip

\begin{proof}
[\textbf{Proof of Theorem \ref{th2.1}}]In our proof we first show that
$\sigma\left(  A_{i}\right)  =\sigma$ for every component of $A_{i}$ and that
conditions (\ref{eqr}) and (\ref{eqc}) hold for each component of $A$. Let
$A_{1},\ldots,A_{k}$ be the components of $A.$ For each $i\in\left[  k\right]
,$ by inequality (\ref{mainin}) we have
\[
\sigma^{2r}\left(  A_{i}\right)  \left\vert R\left(  A_{i}\right)  \right\vert
\geq w_{A_{i}}^{2r+1}\left(  R\left(  A_{i}\right)  \right)  ,
\]
and so%
\[
\sigma^{2r}\left\vert R\right\vert \geq%
{\textstyle\sum\limits_{i\in\left[  k\right]  }}
\sigma^{2r}\left(  A_{i}\right)  \left\vert R\left(  A_{i}\right)  \right\vert
\geq%
{\textstyle\sum\limits_{i\in\left[  k\right]  }}
w_{A_{i}}^{2r+1}\left(  R\left(  A_{i}\right)  \right)  =w_{A}^{2r+1}\left(
R\right)  .
\]
Therefore, condition (\ref{eqr}) implies that $\sigma\left(  A_{i}\right)
=\sigma$ for all $i\in\left[  k\right]  .$ We see also that condition
(\ref{eqr}), and likewise condition (\ref{eqc}), holds for every component of
$A;$ hence, we can assume that $A$ is the sole component, i.e., $A$ is
connected. To finish the proof, we have to show that $A$ is regular.

Since $\sigma^{2}w_{A}^{2s-1}\left(  R\right)  \geq w_{A}^{2s+1}\left(
R\right)  >0$ for all $s\geq1,$ condition (\ref{eqr}) implies that
\[
\sigma^{2}\left(  A\right)  \left\vert R\right\vert =w_{A}^{3}\left(
R\right)  =%
{\textstyle\sum\limits_{i\in R}}
w_{A}^{3}\left(  i\right)  =%
{\textstyle\sum\limits_{i,k\in R,j\in C}}
a_{ij}a_{kj}%
\]
Hence, $\mathbf{j}_{m}$ is an eigenvector of $AA^{\ast}$ to $\sigma^{2}$ and
so, for each $i\in R,$%
\[
\sigma^{2}=%
{\textstyle\sum\limits_{j\in C,k\in R}}
a_{ij}a_{kj}=%
{\textstyle\sum\limits_{j\in C}}
a_{ij}c_{j}.
\]
Let
\[
\delta_{R}=\min\limits_{i\in R}r_{i},\text{ \ \ }\delta_{C}=\min\limits_{i\in
C}c_{i},\text{ \ \ }\Delta_{R}=\max\limits_{i\in R}r_{i},\text{ \ \ }%
\Delta_{C}=\max\limits_{i\in C}c_{i},
\]
and select $s\in R$ such that $r_{s}=\delta_{R}$. Then
\begin{equation}
\sigma^{2}=%
{\textstyle\sum\limits_{j\in C,k\in R}}
a_{sj}a_{kj}=%
{\textstyle\sum\limits_{j\in C}}
a_{sj}c_{j}\leq\Delta_{C}%
{\textstyle\sum\limits_{j\in C}}
a_{sj}=\Delta_{C}\delta_{R}. \label{eqs}%
\end{equation}
Likewise, we see that $\sigma^{2}\geq\delta_{C}\Delta_{R}.$ Applying the same
argument to $A^{\ast}$ we find that
\[
\Delta_{C}\delta_{R}\leq\sigma^{2}\leq\Delta_{R}\delta_{C}.
\]
Therefore, $\Delta_{C}\delta_{R}=\sigma^{2}=\Delta_{R}\delta_{C}.$ This
implies that equality holds in (\ref{eqs}), and so $c_{j}=\Delta_{C}$ whenever
$a_{sj}\neq0.$ Likewise, we see that if $t\in R$ is such that $c_{t}%
=\Delta_{C},$ then $r_{j}=\delta_{R}$ whenever $a_{jt}\neq0.$ Since $A$ is
connected $r_{1}=\cdots=r_{m}$ and $c_{1}=\cdots=c_{n},$ completing the proof.
\end{proof}

\bigskip

\begin{proof}
[\textbf{Proof of Theorem \ref{th3}}]The implications \emph{(i)}$\implies
$\emph{(ii)}$\implies$\emph{(iii) }are obvious, so we shall focus on
\emph{(iii)}$\implies$\emph{(i)}$.$ As in the proof of Theorem \ref{th2.1}, we
first reduce the argument to connected matrices. Let $A_{1},\ldots,A_{k}$ be
the components of $A.$ By Theorem \ref{th1}, for each $s\in\left[  k\right]
,$ we have
\[
\sigma\left(  A_{s}\right)  \sqrt{w_{A_{s}}^{r}\left(  R\left(  A_{s}\right)
\right)  w_{A_{s}}^{r}\left(  C\left(  A_{s}\right)  \right)  }\geq%
{\textstyle\sum\limits_{i\in R\left(  A_{s}\right)  ,j\in C\left(
A_{s}\right)  }}
a_{ij}\sqrt{w_{A_{s}}^{r}\left(  i\right)  w_{A_{s}}^{r}\left(  j\right)  },
\]
and, using the Cauchy-Schwarz inequality,%
\begin{align*}
\sigma\sqrt{w_{A}^{r}\left(  R\right)  w_{A}^{r}\left(  C\right)  }  &
=\sigma\sqrt{%
{\textstyle\sum\limits_{s\in\left[  k\right]  }}
w_{A_{s}}^{r}\left(  R\left(  A_{s}\right)  \right)
{\textstyle\sum\limits_{s\in\left[  k\right]  }}
w_{A_{s}}^{r}\left(  C\left(  A_{s}\right)  \right)  }\\
&  \geq\sigma%
{\textstyle\sum\limits_{s\in\left[  k\right]  }}
\sqrt{w_{A_{s}}^{r}\left(  R\left(  A_{s}\right)  \right)  w_{A_{s}}%
^{r}\left(  C\left(  A_{s}\right)  \right)  }\\
&  \geq%
{\textstyle\sum\limits_{s\in\left[  k\right]  }}
\sigma\left(  A_{s}\right)  \sqrt{w_{A_{s}}^{r}\left(  R\left(  A_{s}\right)
\right)  w_{A_{s}}^{r}\left(  C\left(  A_{s}\right)  \right)  }\\
&  \geq%
{\textstyle\sum\limits_{s\in\left[  k\right]  }}
{\textstyle\sum\limits_{i\in R\left(  A_{s}\right)  ,j\in C\left(
A_{s}\right)  }}
a_{ij}\sqrt{w_{A_{s}}^{r}\left(  i\right)  w_{A_{s}}^{r}\left(  j\right)  }\\
&  =%
{\textstyle\sum\limits_{i\in R,j\in C}}
a_{ij}\sqrt{w_{A}^{r}\left(  i\right)  w_{A}^{r}\left(  j\right)  }.
\end{align*}

Therefore, condition (\ref{eq}) implies that $\sigma\left(  A_{i}\right)
=\sigma$ for all $i\in\left[  k\right]  .$ We see also that condition
(\ref{eq}) holds for every component of $A;$ hence, we can assume that $A$ is
the sole component, i.e., $A$ is connected. To finish the proof we must show
that $A$ is regular.

Let $B=\left(  b_{ij}\right)  $ be defined as a block matrix
\[
B=\left(
\begin{array}
[c]{cc}%
0 & A\\
A^{\ast} & 0
\end{array}
\right)  .
\]
It is known (\cite{HoJo88}, p. 418) that the positive eigenvalues of $B$ are
the nonzero singular values of $A.$ Set for convenience $R=\left[  m\right]  $
and $C=\left[  m+1..m+n\right]  .$ By induction on $r$ it is easy to see that
for every $r\geq0$ and for each $i\in\left[  m+n\right]  ,$ the value
$w_{A}^{r+1}\left(  i\right)  $ is equal to the $i$th row sum of $B^{r}.$

Let%
\begin{align*}
x_{i}  &  =\sqrt{w_{\left\vert A\right\vert }^{r}\left(  i\right)  }\text{ for
}i\in\left[  m\right]  ,\text{ \ \ }y_{i}=\sqrt{w_{\left\vert A\right\vert
}^{r}\left(  i+m\right)  }\text{ for }i\in\left[  n\right]  ,\\
\mathbf{x}  &  =\left(  x_{1},\ldots,x_{m}\right)  ,\text{ \ \ }%
\mathbf{y}=\left(  y_{1},\ldots,y_{n}\right)  ,\text{ \ \ }\mathbf{z}=\left(
x_{1},\ldots,x_{m},y_{1},\ldots,y_{n}\right)  .
\end{align*}
Our main goal is to show that
\begin{equation}
w_{\left\vert A\right\vert }^{r}\left(  1\right)  =\cdots=w_{\left\vert
A\right\vert }^{r}\left(  m\right)  ,\text{ \ \ }w_{\left\vert A\right\vert
}^{r}\left(  m+1\right)  =\cdots=w_{\left\vert A\right\vert }^{r}\left(
m+n\right)  \label{eqw}%
\end{equation}

Equation (\ref{eq}) and the Rayleigh principle imply that $\mathbf{z}$ is an
eigenvector of $B$ to $\sigma;$ hence $\mathbf{z}$ is an eigenvector of
$B^{r}$ to $\sigma^{r}.$ Assume first that $r$ is even, say $r=2k.$ We have
\[
B^{2k-1}\mathbf{z}=\left(
\begin{array}
[c]{cc}%
0 & \left(  AA^{\ast}\right)  ^{k-1}A\\
\left(  A^{\ast}A\right)  ^{k-1}A^{\ast} & 0
\end{array}
\right)  \mathbf{z}=\sigma^{2k-1}\mathbf{z},
\]
and so,
\begin{align}
\sigma^{2k-1}x_{i}  &  =%
{\textstyle\sum\limits_{j\in m+1}^{m+n}}
b_{ij}y_{j-m}\text{ \ \ for }i\in\left[  m\right]  ,\label{eq1}\\
\sigma^{2k-1}y_{i-m}  &  =%
{\textstyle\sum\limits_{j\in\left[  m\right]  }}
b_{ij}x_{j}\text{ \ \ \ for }i\in\left[  m+1..n\right]  . \label{eq2}%
\end{align}
Let%
\[
\delta_{R}=\min\limits_{i\in\left[  m\right]  }w_{A}^{2k}\left(  i\right)
,\text{ \ \ }\delta_{C}=\min\limits_{i\in C}w_{A}^{2k}\left(  i\right)
,\text{ \ \ }\Delta_{R}=\max\limits_{i\in R}w_{A}^{2k}\left(  i\right)
,\text{ \ \ }\Delta_{C}=\max\limits_{i\in C}w_{A}^{2k}\left(  i\right)  .
\]
Select $s\in\left[  m\right]  $ such that $x_{s}=\delta_{R}$. Then, by
(\ref{eq1}),%
\begin{equation}
\sigma^{2k-1}\sqrt{\delta_{R}}=\sigma^{2k-1}x_{s}=%
{\textstyle\sum\limits_{j=m+1}^{n}}
b_{sj}y_{j-m}\leq\sqrt{\Delta_{C}}%
{\textstyle\sum\limits_{j=m+1}^{n}}
b_{sj}=\sqrt{\Delta_{C}}\delta_{R}, \label{eqs1}%
\end{equation}
and so $\sigma^{2k-1}\leq\sqrt{\Delta_{C}\delta_{R}}.$ Likewise, we see that
$\sigma^{2k-1}\geq\sqrt{\Delta_{R}\delta_{C}}.$ Applying the same argument to
equation (\ref{eq2}), we find that
\[
\sqrt{\Delta_{C}\delta_{R}}\leq\sigma^{2k-1}\leq\sqrt{\Delta_{R}\delta_{C}}.
\]
Therefore, $\Delta_{C}\delta_{R}=\sigma^{2k-1}=\Delta_{R}\delta_{C}.$ This
implies that equality holds in (\ref{eqs1}), and so $w_{A}^{2k}\left(
j\right)  =\Delta_{C}$ whenever $b_{sj}\neq0.$ Likewise, we see that if
$t\in\left[  m+1..n\right]  $ is such that $w_{A}^{2k}\left(  i\right)
=\Delta_{C},$ then $w_{A}^{2k}\left(  j\right)  =\delta_{R}$ whenever
$b_{jt}\neq0.$ Since $A$ is connected, (\ref{eqw}) holds for even $r$. The
proof of (\ref{eqw}) for odd $r\mathbf{\ goes}$ along the same lines and we
omit it.

Note that (\ref{eqw}) implies that $\mathbf{j}_{m+n}$ is an eigenvector to $B$
and thus all row and column sums of $A$ are equal, completing the proof.
\end{proof}

\bigskip

\begin{proof}
[\textbf{Proof of Theorem \ref{th4}}]Suppose $A$ is regular. Schur's
inequality \cite{Sch11}
\[
\sigma^{2}\left(  A\right)  \leq\max_{i\in R,j\in C}r_{i}c_{j},
\]
implies that $\sigma\left(  A\right)  \leq\Sigma\left(  A\right)  /\sqrt{nm}.$
In view of (\ref{secin1}), we deduce that $\sigma\left(  A\right)
=\Sigma\left(  A\right)  /\sqrt{nm}.$

Suppose now $\sigma\left(  A\right)  =\Sigma\left(  A\right)  /\sqrt{nm}.$ We
have%
\[
\sigma^{2}\left(  A\right)  m\geq\left\langle AA^{\ast}\mathbf{j}%
_{m},\mathbf{j}_{m}\right\rangle =%
{\textstyle\sum\limits_{i,k\in R,j\in C}}
a_{ij}a_{kj}=%
{\textstyle\sum\limits_{j\in C}}
c_{j}^{2}\geq\frac{1}{n}\left(  \Sigma\left(  A\right)  \right)  ^{2}%
=\sigma^{2}\left(  A\right)  m,
\]
implying that $c_{1}=\cdots=c_{n}$. Likewise we find that $r_{1}=\cdots
=r_{m},$ completing the proof.
\end{proof}

\bigskip

\textbf{Concluding remark}

It seems a challenging problem to investigate the cases of equality in () and
() for arbitrary matrices.

\end{document}